\pdfoutput=1
\RequirePackage{ifpdf}
\ifpdf % We are running pdfTeX in pdf mode
\documentclass[pdftex]{sigma}
\else
\documentclass{sigma}
\fi

\numberwithin{equation}{section}

\newtheorem{Theorem}{Theorem}[section]
\newtheorem{Proposition}[Theorem]{Proposition}

\begin{document}

\allowdisplaybreaks

\newcommand{\arXivNumber}{1705.07625}

\renewcommand{\thefootnote}{}

\renewcommand{\PaperNumber}{088}

\FirstPageHeading

\ShortArticleName{Variations for Some Painlev\'e Equations}

\ArticleName{Variations for Some Painlev\'e Equations\footnote{This paper is a~contribution to the Special Issue on Algebraic Methods in Dynamical Systems. The full collection is available at \href{https://www.emis.de/journals/SIGMA/AMDS2018.html}{https://www.emis.de/journals/SIGMA/AMDS2018.html}}}

\Author{Primitivo B.~ACOSTA-HUM\'ANEZ~$^{\dag\ddag}$, Marius~VAN DER PUT~$^\S$ and Jaap~TOP~$^\S$}

\AuthorNameForHeading{P.B.~Acosta-Hum\'anez, M.~van der Put and J.~Top}

\Address{$^\dag$~School of Basic and Biomedical Sciences, Universidad Sim\'on Bol\'{\i}var,\\
\hphantom{$^\dag$}~Barranquilla, Colombia}
\Address{$^\ddag$~Instituto Superior de Formaci\'on Docente Salom\'e Ure\~na - ISFODOSU,\\
\hphantom{$^\ddag$}~Santiago de los Caballeros, Dominican Republic}
\EmailD{\href{mailto:primitivo.acosta-humanez@isfodosu.edu.do}{primitivo.acosta-humanez@isfodosu.edu.do}}

\Address{$^\S$~Bernoulli Institute, University of Groningen, Groningen, The Netherlands}
\EmailD{\href{mailto:m.van.der.put@rug.nl}{m.van.der.put@rug.nl}, \href{mailto:j.top@rug.nl}{j.top@rug.nl}}

\ArticleDates{Received November 01, 2018, in final form November 05, 2019; Published online November 09, 2019}

\Abstract{This paper first discusses irreducibility of a Painlev\'e equation $P$. We explain how the Painlev\'e property is helpful for the computation of special classical and algebraic solutions. As in a paper of Morales-Ruiz we associate an autonomous Hamiltonian $\mathbb{H}$ to a~Painlev\'e equation~$P$. Complete integrability of $\mathbb{H}$ is shown to imply that all solutions to~$P$ are classical (which includes algebraic), so in particular~$P$ is solvable by ``quadratures''. Next, we show that the variational equation of~$P$ at a given algebraic solution coincides with the normal variational equation of~$\mathbb{H}$ at the corresponding solution. Finally, we test the Morales-Ramis theorem in all cases $P_{2}$ to $P_{5}$ where algebraic solutions are present, by showing how our results lead to a quick computation of the component of the identity of the differential Galois group for the first two variational equations. As expected there are no cases where this group is commutative.}


\Classification{33E17; 34M55}

\renewcommand{\thefootnote}{\arabic{footnote}}
\setcounter{footnote}{0}

\vspace{-3mm}

\section{Introduction and summary}
The interesting idea of J.-A.~Weil to apply the Morales-Ramis theorem to Painlev\'e equations was initiated in~\cite{Mo2}. It is also the subject of more recent papers \cite{A,HS,St3,St2,St4,SC}. The Hamiltonian $H$ of a~Painlev\'e equation $x''=R(x',x,t)$ depends on `the time' $t$. In order to apply the Morales-Ramis theorem, $H$ is changed into a time-independent Hamiltonian $\mathbb{H}=H+e$. Our first main result (Proposition~\ref{2.1}) states that complete integrability for $\mathbb{H}$ implies that all solutions of the equation $x''=R(x',x,t)$ are classical solutions in the sense of H.~Umemura \cite{U3,U1,U2,UW} (this includes algebraic functions). In fact, one may state that such an equation $x''=R(x',x,t)$ is {\it not} considered as a true Painlev\'e equation. This is in agreement with \cite{ZF}, see also Section~\ref{Sect2} below.

The second main result of this paper (Proposition~\ref{prop3.1}) claims that the normal variational equation(s) of $\mathbb{H}$ along a given explicit solution are shown to be equivalent to the variational equation(s) for $x''=R(x',x,t)$ along a given solution. To fix notations, we recall notations for second-order Painlev\'{e} equations here, taken from~\cite{OO2} and~\cite{APT}:
\begin{alignat*}{3}
& P_1\colon \quad && y''= 6y^2+t,& \\
& P_2(\alpha)\colon \quad && y''= 2y^3+ty+\alpha, & \\
& P'_3(\alpha,\beta,\gamma,\delta)\colon \quad && y''=y'^2/y -y'/t+ \big(\alpha y^2+\gamma y^3\big)/\big(4t^2\big)+ \beta/(4t)+\delta/(4y), &\\
& P_4(\alpha,\beta)\colon \quad && y''=y'^2/y +3y^3/2+4ty^2+2\big(t^2-\alpha\big)y +\beta/y, &\\
& P_5(\alpha,\beta,\gamma,\delta)\colon \quad && y''=\left(\frac{1}{2y} +\frac{1}{y-1}\right)y'^2-y'/t+(y-1)^2(\alpha y+\beta/y)/\big(t^2\big)\\
&&& \hphantom{y''=}{} +\gamma y/t+\delta y(y+1)/(y-1).&
\end{alignat*}
We will not use a formula for $P_6$. In the case of $P_3'$ there is a refinement:
\begin{gather*}
P'_3(D_6)=P'_3(\alpha,\beta,\gamma,\delta)\qquad \text{with}\quad \gamma\delta\neq 0,\\
P'_3(D_7)=P'_3(\alpha,\beta,\gamma,\delta)\qquad \text{with}\quad (\delta=0, \beta\neq 0)\quad \text{or}\quad (\gamma=0,\alpha\neq 0),\\
P'_3(D_8)=P'_3(\alpha,\beta,0,0)\qquad \text{with}\quad \alpha\beta\neq 0.
\end{gather*}
Finally we mention
\[
\deg P_5(\theta_0,\theta_1)\colon \
y''=\frac{1}{2}\left(\frac{1}{y}+\frac{1}{y-1}\right)(y')^2-\frac{y'}{t}+\frac{2(y-1)\theta_0^2}{yt^2}-\frac{2y\theta_1^2}{(y-1)t^2}+8y(y-1).
\]
In \cite[Section~3]{APT} it is explained how this relates to $P_5\big(\theta_1^2/2,-\theta_0^2/2,-2,0\big)$.

For the equations $P_{2}$ to $P_5$, there is a convenient list in~\cite{OO2} of all cases with algebraic solutions, up to B\"acklund transformations. The Hamiltonians arising from the equations in this list are not completely integrable, as follows from Propositions~\ref{2.1} and \ref{prop3.1}. One expects, in accordance with \cite{MR, MRS}, that the variational equations in such cases produce differential Galois groups $G$ such that $G^o$, the component of the identity, is not abelian. We verify this explicitly for all items in the list of~\cite{OO2}. The results and some comments regarding them are (see Section~\ref{Sect4}):
\begin{itemize}\itemsep=0pt
\item[(1)] The first variational equation produces, for almost all cases, $G={\rm SL}_2$.
\item[(2)] In some cases the first variational equation produces the differential Galois group with component of the identity~$\mathbb{G}_m$. The second variational equation produces an extension of this group by a unipotent group~$\mathbb{G}_a^m$. The action by conjugation of $\mathbb{G}_m$ on~$\mathbb{G}_a^m$ is not trivial. Hence $G^o$ is not abelian.

\item[(3)] We elaborate the interesting case Section~\ref{4.7} which discusses $P_5(a,-a,0,\delta)$ and $y=-1$ with first variational equation $v''=-t^{-1}v' +\big(8at^{-2}+\frac{1}{2}\delta \big)v$. It follows from the monodromy theorem~\cite[Proposition~8.12(2)]{vdP-Si} that its differential Galois group is ${\rm SL}_2$ unless $a=\frac{(2n+1)^2}{32}$ with $n\in \mathbb{Z}$. For these special values of~$a$ the differential Galois group is~$\mathbb{G}_m$. Again, for these special cases, the second variational produces a $G$ with non abelian component of the identity.
In \cite{St2} a $P_5$ equation with different parameters (but equivalent by B\"acklund transformations) is studied and the same special values are found. Our methods (specifically, the use of the monodromy theorem) simplify, compared to earlier similar results by Stoyanova et al.~\cite{St3,St2,St4, SC} the determination of the Galois group associated to such a~variational equation.
The special values for $a$ can be explained as follows. There is a standard isomonodromy family corresponding to $P_5$, see \cite{JM}. Let $\pm \frac{\theta_0}{2}$ and $\pm \frac{\theta_1}{2}$ denote the local exponents of this family for the regular singularities~$0$ and~$1$. Then $\theta_0-\theta_1=\sqrt{8a}$. Thus the special values for~$a$ correspond to a type of resonance between the regular singularities at~0 and at~1.

\item[(4)] (This observation is in part inspired by a discussion of one of us with Juan J.~Morales-Ruiz; we thank him for his question.) Each of the Painlev\'{e} equations is induced by isomonodromy of some family of order~$2$ linear differential equations~\cite{vdP-Sa}. The possible singular points of such a family are $0$, $1$, $\infty$ with prescribed singularity. The first variational equation happens to have the same type of possible singularities; this observation is made in various examples discussed in Section~\ref{Sect4}.
\end{itemize}

R.~Fuchs' problem, see \cite{OO2,St1}, also concerns algebraic solutions of Painlev\'e equations. The second-order linear differential equations resulting from this problem seem to be unrelated to the first variational equations.

In Section~\ref{one} we observe that if a second-order equation $R$ has the Painlev\'{e} property and moreover is reducible, then the induced first-order differential equation~$Q$ has the Painlev\'e pro\-per\-ty, too. The classification of first-order equations with the Painlev\'e pro\-perty has consequences for the special solutions of $R$, as will be explained on p.~\pageref{pagefive}.

\section{Reducibility and special solutions}\label{one}
Consider a Painlev\'e equation $x''=R(x',x,t)$ with fixed parameters. Let ${\rm den}$ be the denominator of $R(x',x,t)$ seen as element of the field of fractions of $\mathbb{C}(t)[x',x]$. Then $\mathcal{D}:=\mathbb{C}(t)\big[x',x,\frac{1}{{\rm den}}\big]$ is a differential algebra with respect to the differentiation given by $t'=1$, $x'=x'$, $x''=R(x',x,t)$.

A solution of a Painlev\'{e} equation is called ``known'' or ``reducible'' if it is obtained from solutions of linear equations, first-order equations, and abelian integrals.
Many results on solutions of Painlev\'{e} equations are known, for example due to the Belarusian school~\cite{G}. A~Painlev\'{e} equation is called reducible if it has a reducible solution; otherwise it is called irreducible. We note that a different definition of reducible second-order equation appears in~\cite{CW}. A~deep result of the Japanese school translates the non-existence of reducible solutions (so, irreducibility of the equation) into
\begin{itemize}\itemsep=0pt
\item[(1)] {\it there are no algebraic solutions}, and
\item[(2)] {\it For every differential field extension $K\supset \mathbb{C}(t)$ the ring $K\otimes\mathcal{D}$ has no principal differential ideal $\neq (0),(1)$}.
\end{itemize}

For this subject we refer to \cite{OKSO,U3, U1,U2} and in particular to~\cite[Appendix~A]{NP}. In fact it is known that condition (2) can be replaced by the simpler condition:
\begin{itemize}\itemsep=0pt
\item[$(2')$] {\it $\mathcal{D}$ has no principal differential ideal $\neq (0),(1)$.}
\end{itemize}

We now discuss how to verify these two conditions.

Concerning (1): By the Painlev\'e property, an algebraic solution can only be ramified above the fixed singularities. These are $t=\infty$ in the cases $P_1$, $P_2,$ and $P_4$. Hence here an algebraic solution must be rational. It is easily seen that no solution of $P_1$ in $\mathbb{C}(t)$ exists. For $P_3$ and $P_5$ an algebraic solution can only ramify above $t=0,\infty$. The equations~$P_6$ have many algebraic solutions; they ramify above $t=0,1,\infty$.

Concerning $(2')$: The algebra $\mathcal{D}$ has unique factorization and one easily verifies that every prime factor $Q$ of $F$ such that $(F)$ is a differential ideal, generates again a~differential ideal. Thus the equation is reducible if and only if $\mathcal{D}$ has a prime ideal $(Q)$ of height one which is invariant under differentiation.

\label{pagefive}Now $Q(x',x,t)=0$ is a first-order differential equation. It is well known that $x''=R(x',x,t)$ has the Painlev\'e property. The solutions of $Q(x',x,t)=0$ are also solutions of this Painlev\'e equation. Therefore $Q(x',x,t)=0$ itself has the Painlev\'e property. A classical result (see \cite{M,Mun-vdP,NNPT} for modern proofs and references to some of the rich classical literature) implies that $Q(x',x,t)$ has one of the following properties:
\begin{itemize}\itemsep=0pt
\item[(i)] Genus 0; it is a Riccati equation; thus $x'=a+bx+cx^2$ with $a,b,c\in \mathbb{C}(t)$.
\item[(ii)] Genus 1; it is a Weierstrass equation; thus it is equivalent to $(x')^2=f\cdot \big(x^3+ax+b\big)$, where $f\neq 0$ is algebraic over $\mathbb{C}(t)$, and $a,b\in \mathbb{C}$ are such that the equation $y^2=x^3+ax+b$ represents an elliptic curve.
\item[(iii)] Genus $>1$; after a finite extension of $\mathbb{C}(t)$, the equation is equivalent to the equation $x'=0$. This is equivalent to the statement: all solutions of $Q(x',x,t)=0$ are contained in a fixed finite extension of $\mathbb{C}(t)$.
\end{itemize}

We will call the above three cases ``special classical solutions'' of the Painlev\'e equation. {\it We conclude that the Painlev\'e equation is irreducible if it has no special classical solutions}.

We observe that all of the cases of special solutions have genus~0. The corresponding second-order linear differential equation has at most singularities where the Painlev\'e equation has fixed singularities (for~$P_6$ these are the points $0$, $1$, $\infty$; for $P_5$, $\deg P_5$, $P_3$ the points $0$, $\infty$; for $P_4$, $P_{2,{\rm FN}}=P_{34}$, $P_2$, $P_1$ the point~$\infty$).

\section[Complete integrability for the Hamiltonian $\mathbb{H}$]{Complete integrability for the Hamiltonian $\boldsymbol{\mathbb{H}}$}\label{Sect2}

$x''=R(x',x,t)$ is again some Painlev\'e equation with fixed parameters. There is a Hamiltonian function $H(y,x,t)$ related to the given Painlev\'e equation. There are various possibilities for $H$ but we assume that it is a rational in the variables $y$, $x$, $t$ and moreover is polynomial of degree~2 in the variable~$y$. The usual equations are:
\[x'(t)=\frac{\partial H}{\partial y}(y(t),x(t),t) \qquad \mbox{and} \qquad y'(t)=-\frac{\partial H}{\partial x}(y(t),x(t),t). \]
Since $H$ is a polynomial of degree two in $y$, the first equation can be used to write $y(t)$ as a~rational expression in $x'(t)$, $x(t)$ and $t$. Substitution of this expression for $y(t)$ in the second equation will produce an explicit second-order equation for $x(t)$ and this is the given one $x''=R(x',x,t)$.

Now $H$ depends on the time $t$. One wants to apply the Morales-Ramis theorem concerning complete integrability. This leads to a choice of a new Hamiltonian $\mathbb{H}(y,x,z,e)=H(y,x,z)+e$ which depends on two pairs of variables $y$, $x$ and $z$, $e$. The new equations are
\[ x'(t)=\frac{\partial \mathbb{H}}{\partial y}, \qquad y'(t)=-\frac{\partial \mathbb{H}}{\partial x}, \qquad z'(t)=\frac{\partial \mathbb{H}}{\partial e}=1,\qquad e'(t)=-\frac{\partial \mathbb{H}}{\partial z}.\]

\begin{Proposition}\label{2.1} Suppose that $\mathbb{H}$ is completely integrable. Then all solutions of $x''=R(x',x,t)$ are classical $($including algebraic$)$. In particular the equation is reducible.
\end{Proposition}
\begin{proof} $\mathbb{H}$ is a first integral. There is an independent first integral $E(y,x,z,e)$. We suppose that $E$ is a rational (or algebraic) function of the 4 variables. Now we replace the $e$ in $E$ by $-H(y,x,z)$. The result is a first integral $F(y,x,z)$ for $H$ and a rational (or algebraic) function $G=G(x',x,t)$ such that $G(x'(t),x(t),t)$ is independent of $t$ for every solution $x(t)$ of the Painlev\'e equation $x''=R(x',x,t)$. Then, taking the derivative with respect to $t$, one finds that the expression
$R(x',x,t)\frac{\partial G}{\partial x'}+x'\frac{\partial G}{\partial x}+\frac{\partial G}{\partial t}$ is zero on every solution $(x'(t),x(t),t)$ of the Painlev\'e equation. It follows that this expression itself is zero.

Consider, as before, the differential algebra $\mathcal{D}:=\mathbb{C}(t)\big[x',x,\frac{1}{{\rm den}}\big]$ with derivation $F\mapsto F'$ given by $t'=1$, $(x)'=x'$, $(x')'=R(x',x,t)$. Assume (for convenience) that $G$ is rational. Thus $G$ lies in the field of fractions $Qt(\mathcal{D})$ of $\mathcal{D}$ and $G'=0$. Let $L\subset Qt(\mathcal{D})$ denote the field of constants. Then $L(x',x,t)$ equals $Qt(\mathcal{D})$ and the transcendance degree of $L(x',x,t)\supset L$ is 2, because $G\not \in \mathbb{C}$. Therefore there is an irreducible polynomial $Q\in L(t)[S,T]$ such that $Q(x',x)=0$. The coefficients of $Q$ lie in $\mathcal{D}[\frac{1}{U}]$ for a suitable element~$U$.

The solutions of the Painlev\'e equation correspond to $\mathbb{C}(t)$-linear differential homomorphism $\phi\colon\mathcal{D}\rightarrow {\rm Mer}$, where ${\rm Mer}$ denotes the differential field of the multivalued meromorphic functions on, say, $\mathbb{C}\setminus \{0,1\}$. Indeed, the homomorphism $\phi$ corresponds to the solution $\phi (x)\in {\rm Mer}$.

If $\phi(U)\neq 0$, then $\phi(Q)(S,T)$ makes sense. Since the coefficients of $Q$ are rational functions in~$t$ with `constant' coefficients, one has $\phi(Q)(S,T)\in \mathbb{C}(t)[S,T]$. Moreover $\phi(Q)(\phi(x)',\phi(x))=0$, which means that the solution $\phi(x)\in {\rm Mer}$ satisfies a~first-order differential equation, which has again the Painlev\'e property. We conclude that $\phi(x)$ is a~special classical solution or an algebraic solution.

Consider finally a $\phi$ such that $\phi(U)=0$. Then $\phi$ is also zero on a prime differential ideal of~$\mathcal{D}$ containing~$U$. If this is a principal ideal, then $\phi(x)$ is a~special classical solution. If this is a~maximal ideal, then the solution $\phi(x)$ is algebraic. \end{proof}

We note that Proposition~\ref{2.1} is in agreement with a main result of~\cite{ZF}: {\it The Hamilton system $\mathbb{H}$, associated with any of the equations $P_1$--$P_{6}$, does not admit any first integral which is an algebraic function of $x$, $y$, $z$, $e$ and independent of $\mathbb{H}$, except in the following cases: {\rm (a)}~$\alpha=\gamma=0$ in $P_{3}$, {\rm (b)}~$\beta=\delta=0$ in $P_{3}$ and {\rm (c)}~$\gamma=\delta=0$ in $P_5$.}

It is well known that in the cases (a)--(c) all solutions are obtained by ``quadratures''. The `first integrals' (this means here an $F$ such that $(F)\subset \mathcal{D}$ is a (prime) differential ideal) are actually known. Namely for $P_{3}$ with $\beta=\delta=0$ they are $t^2 (x')^2+2txx'-\big(C+2\alpha tx+\gamma t^2x^2\big)x^2$ (with arbitrary $C$ and a similar formula for the case $\alpha=\gamma=0$); see slide~47 of the 2002 lecture~\cite{C2} by Clarkson and also the Russian paper~\cite{Luk67}. For $P_5$ with $\gamma =\delta=0$ the `first integrals' are $t^2(x')^2-(x-1)^2\big(2\alpha x^2+Cx-2\beta\big)$. This is, e.g., stated on slide~48 of~\cite{C2}, see also~\cite{Luk68}.

One observes that the above `first integrals' are order one differential equations having the Painlev\'e property. They have genus~$0$.

\section{Several variational equations}\label{section3}
Suppose that an algebraic solution $x_0=x_0(t)$ of $x''=R(x',x,t)$ is given. The variational equation VEP for the Painlev\'e equation is given by the following formalism. Put $x=x_0(t)+\epsilon v$ with $\epsilon ^2=0$. Substitution yields the equation $x''_0+\epsilon v''=R(x_0'+\epsilon v',x_0+\epsilon v,t)$. The coefficient of $\epsilon$ in this equation is a second-order linear equation for $v$. Explicitly, VEP is the equation
\[ v''=\frac{\partial R}{\partial x'}(x'_0,x_0,t)v'+\frac{\partial R}{\partial x}(x_0',x_0,t)v. \]

The algebraic solution $x_0$ produces an algebraic solution $y_0=y_0(t)$, $x_0=x_0(t)$ for the Hamiltonian equations for~$H$. The variational equation VEH for this Hamiltonian equation is defined by the following formalism. Put $y=y_0+\epsilon w$, $x=x_0+\epsilon v$ with $\epsilon ^2=0$ in the two Hamilton equations. Thus
\[x_0'+\epsilon v' =\frac{\partial H}{\partial y}(y_0+\epsilon w,x_0+\epsilon v,t) \qquad \text{and} \qquad y_0'+\epsilon w'=- \frac{\partial H}{\partial x}(y_0+\epsilon w,x_0+\epsilon v,t). \]

The coefficients of $\epsilon$ in these equations yield linear differential equations for $w$ and $v$ of order one. Moreover, the first equation can be used to eliminate $w$ as
linear expression in $v$ and $v'$. Thus we obtain a second-order homogenous differential equation for $v$ which coincides of course with the earlier VEP.

An algebraic solution $x_0=x_0(t)$ for the Painlev\'e equation yields for $\mathbb{H}$ the solution $y_0=y_0(t)$, $x_0=x_0(t)$, $z_0(t)=t$, $e_0=e_0(t)=-\int \frac{\partial \mathbb{H}}{\partial z}{\rm d}t$. For the algebraic solutions of $P_2,\ldots, P_5$ and the Hamiltonians as given in~\cite{OO2}, the function~$e_0$ turns out to be algebraic. We have not verified this for the case of $P_6$. If in such a situation a transcendental $e_0$ occurs, then the Morales-Ramis theory is still valid, and the VE$\mathbb{H}$ and NVE$\mathbb{H}$ still make sense. However, in that case these equations are considered over a differential field which is larger than the algebraic closure of~$\mathbb{C}(t)$.

The variational equation VE$\mathbb{H}$ for $\mathbb{H}$ along this solution is obtained by the following formalism. Put $y=y_0+\epsilon w$, $x=x_0+\epsilon v$, $z_0=t+\epsilon a$, $e=e_0+\epsilon b$ with $\epsilon ^2=0$. Substitution of these data into the Hamiltonian equations for~$\mathbb{H}$ yields a systems of rank~4 of linear differential equations of first order. In more detail
\begin{gather*} x_0'+\epsilon v'=\frac{\partial \mathbb{H}}{\partial
 y}(y_0+\epsilon w, \dots, e_0+\epsilon b),\qquad y_0'+\epsilon w'=-\frac{\partial \mathbb{H}}{\partial x}(y_0+\epsilon w,\dots ,e_0+\epsilon b),\\
t'+\epsilon a'=\frac{\partial \mathbb{H}}{\partial e}=1,\qquad e_0'+\epsilon b'=-\frac{\partial \mathbb{H}}{\partial z}(y_0+\epsilon
w,\dots ,e_0+\epsilon b).
\end{gather*}

The normal variational equation NVE$\mathbb{H}$ for $\mathbb{H}$ along
this solution is obtained by taking the three dimensional space
perpendicular to the equation for $a$ and dividing out by the tangent
line of the curve (i.e., the given solution). This means that we
are reduced to the case $a=b=0$.
Substitution of $a=b=0$ in the VE$\mathbb{H}$
yields the VEH, which is equivalent to the VEP.
We note that the VE$\mathbb{H}$ contains the term $e_0$, but NVE$\mathbb{H}$ does not.
This computation proves the following result.

\begin{Proposition}\label{prop3.1} The normal variational equation {\rm NVE}$\mathbb{H}$ of $\mathbb{H}$
coincides with the variational VEP of the Painlev\'e equation $x''=R(x',x,t)$.
\end{Proposition}

\section[VEP for the algebraic solutions of $P_{2},\dots,P_5$]{VEP for the algebraic solutions of $\boldsymbol{P_{2},\dots,P_5}$} \label{Sect4}
We adopt here the list of special solutions of \cite[Theorem~2.1]{OO2}. Further we will, as in that paper, replace the classical $P_{3}$ by $P'_3$. Finally for the degenerate fifth Painlev\'e equation we will use the $\deg P_5$ of our paper~\cite{APT}.

\subsection[$P_{2}(\alpha =0)$ with solution $y=0$]{$\boldsymbol{P_{2}(\alpha =0)}$ with solution $\boldsymbol{y=0}$} \label{section4.1}

$P_{2}$ reads $y''=2y^3+ty+\alpha$. The VEP for $\alpha=0$ and $y=0$ reads $v''=tv$. This is the Airy equation with differential Galois group ${\rm SL}_2$. This is also present in~\cite{A,Mo2,SC}.

We skip the Flaschka--Newell $P_{34}=P_{2,{\rm FN}}$ since it is equivalent to $P_2$.

\subsection[$P_4\big(0,-\frac{2}{9}\big)$ with $y=-\frac{2}{3}t$]{$\boldsymbol{P_4\big(0,-\frac{2}{9}\big)}$ with $\boldsymbol{y=-\frac{2}{3}t}$}\label{4.2}

In this case the VEP reads $v''=t^{-1}v'-\frac{4}{3}t^2v$. A basis of solutions is $\big\{ {\rm e}^{\sqrt{-1/3} t^2},{\rm e}^{-\sqrt{-1/3} t^2}\big\}$ and the differential Galois group is~$\mathbb{G}_m$. The second variational equations (obtained by putting as solution $-\frac{2}{3}t+\epsilon v+\epsilon^2 w$) read
\[ v''-t^{-1}v'+\frac{4}{3}t^2v=0,\qquad w''-t^{-1}w'+\frac{4}{3}t^2w=\frac{3}{2}t^{-1}vv''-\frac{3}{4}t^{-1}(v')^2+3tv^2.\]
For every solution $v_0\neq 0$ of the first equation, the second inhomogeneous equation produces an extension of $\mathbb{G}_m$ by an additive group $\mathbb{G}_a$. For instance, the choice $v_0={\rm e}^{ct^2}$ with $c^2=-\frac{1}{3}$ leads to the equation $w''-t^{-1}w'+\frac{4}{3}t^2w=\big(3ct^{-1}+2t\big){\rm e}^{2ct^2}$. Any special solution $w_0$ involves the ``error function'' $\operatorname{erf}(t)$. An element of $s\in\mathbb{G}_m\cong \mathbb{C}^*$ maps $v_0$ to $sv_0$ and maps $w_0$ to $s^2w_0+$ a~solution of the homogeneous equation. Hence the component of the identity of the differential Galois group of the second variational equation is not commutative. Compare also to the next case, where more details are given.

\subsection[$P_4(0,-2)$ with $y=-2t$]{$\boldsymbol{P_4(0,-2)}$ with $\boldsymbol{y=-2t}$}\label{4.3}

This case is similar to the one discussed in Section~\ref{4.2}. The resulting VEP is $v''-t^{-1}v'-4t^2v=0$. A basis of solutions is $\big\{{\rm e}^{t^2},{\rm e}^{-t^2}\big\}$ and the differential Galois group is $\mathbb{G}_m$. The second VEP is
\[v''-t^{-1}v'-4t^2v=0,\qquad w''-t^{-1}w'-4t^2w=\frac{1}{2} t^{-1} v v''-\frac{1}{4}t^{-1}(v')^2-7tv^2.\]
The Picard--Vessiot field $K$ for this set of equations is an extension of the Picard--Vessiot field $K_0=\mathbb{C}\big(t,{\rm e}^{t^2}\big)\supset \mathbb{C}(t)$ of the first
equation. The extension $K\supset K_0$ is obtained by adding a~particular solution $w_0$ of the second inhomogeneous equation for every solution $v_0=a{\rm e}^{t^2}+b{\rm e}^{-t^2}\neq 0$ of the first equation. The differential Galois group~$G$ of~$K/K_0$ is an unipotent group and can be seen to be $\mathbb{G}_a^3$. The action of $\mathbb{G}_m$ on the solutions $v_0$ induces a~non-trivial action (by conjugation) of $\mathbb{G}_m$ on $\mathbb{G}_a^3$. Indeed, fix a solution $v_0$ of the first equation: $v_0={\rm e}^{t^2}$. Choose $\sigma\in\mathbb{G}_m$ such that $\sigma(v_0)=cv_0$ and $c^2\neq 1$ and a solution $w_0$ of the inhomogeneous equation appearing in the second VEP. Any extension $\tilde{\sigma}$ of $\sigma$ to $G$ satisfies $\tilde{\sigma}(w_0)=c^2w_0+v$ with $v''-t^{-1}v'-4t^2v=0$. This shows that the action of $\mathbb{G}_m$ on $\mathbb{G}_a^3$ is non-trivial. Thus the component of the identity of~$G$ is not commutative.

Stoyanova in~\cite{St4} studies non-integrability of $P_4(1,0)$. After suitable B\"{a}cklund transformations her results agree with our Section~\ref{4.3}.

 \subsection[$P_3'(D_6)(a,-a,4,-4)$ and $y=-t^{1/2}$]{$\boldsymbol{P_3'(D_6)(a,-a,4,-4)}$ and $\boldsymbol{y=-t^{1/2}}$} \label{4.4}
The variational equation reads $v''+\big(\frac{1}{4}t^{-2}+\frac{a}{2}t^{-3/2}-4t^{-1}\big)v=0$. At $t=0$ there is a logarithm present in the local solutions. At $t=\infty$ there is an exponential present in the local solutions. The differential Galois group~$G$, say over $\mathbb{C}\big(t^{1/2}\big)$, contains therefore~$\mathbb{G}_a$ and~$\mathbb{G}_m$. Thus $G$ contains a~Borel subgroup. The operator corresponding to this equation does not factor over~$\mathbb{C}\big(t^{1/2}\big)$. One concludes that $G={\rm SL}_2$.

\subsection[$P_3'(D_7)(0,-2,2,0)$ and $y=t^{1/3}$]{$\boldsymbol{P_3'(D_7)(0,-2,2,0)}$ and $\boldsymbol{y=t^{1/3}}$} \label{4.5}

The variational equation reads $-t^{1/3} v''-\frac{1}{3}t^{-2/3}v'+\big({-}\frac{1}{9}t^{-5/3}+ \frac{3}{2}t^{-1}\big)v=0$. As in Section~\ref{4.4}, there is at $t=0$ a logarithm present and at $t=\infty$ an exponential. The corresponding operator does not factor over the field $\mathbb{C}\big(t^{1/3}\big)$. One concludes that the differential Galois group is ${\rm SL}_2$.

\subsection[$P_3'(D_8)(8h,-8h,0,0)$ and $y=-t^{1/2}$]{$\boldsymbol{P_3'(D_8)(8h,-8h,0,0)}$ and $\boldsymbol{y=-t^{1/2}}$} \label{4.6}

The variational equation is $v''+\big(4ht^{-3/2}+\frac{1}{4}t^{-2}\big)v=0$. Computations similar to those in Sec\-tions~\ref{4.4} and~\ref{4.5} imply that the differential Galois group is ${\rm SL}_2$.

\subsection[$P_5(a,-a,0,\delta)$ and $y=-1$]{$\boldsymbol{P_5(a,-a,0,\delta)}$ and $\boldsymbol{y=-1}$}\label{4.7}

Here the VEP is $v''=-t^{-1}v'+\big(8at^{-2}+\frac{1}{2}\delta\big)v$ and the operator form is $\delta _t^2-\big(8a+\frac{1}{2}\delta t^2\big)$ with $\delta_t:=t\frac{{\rm d}}{{\rm d}t}$. The differential Galois group $G\subset {\rm SL}_2$ depends on~$a$ and~$\delta$.

If $\delta =0$, then this group is a subgroup of $\mathbb{G}_m$. We skip this equation because it is a special case of the degenerate fifth Painlev\'e equation. (Alternatively, observe that it is `quadrature'.)

Suppose that $\delta \neq 0$ (then one usually scales $\delta$ to $-\frac{1}{2}$). We use~\cite{vdP-Si} for some facts and terminology and we use the package DEtools of MAPLE for some computations. The singularity at $t=\infty$ is irregular; its generalized eigenvalues are $\pm \big(\frac{\delta}{2}\big)^{1/2}t$ and the formal monodromy $\gamma$ is~$-{\rm id}$. On a suitable basis of the formal solution space at $t=\infty$, the topological monodromy at $t=\infty$ is the product of $\gamma$ and two Stokes matrices and has the form $\left(\begin{smallmatrix} -1& 0 \\ 0& -1\end{smallmatrix}\right)
\left(\begin{smallmatrix}1& 0\\ e_1& 1\end{smallmatrix}\right) \left(\begin{smallmatrix}1& e_2\\ 0 & 1\end{smallmatrix}\right)$.

Since there is only one other singular point, namely at $t=0$, and since this singularity is regular, the group~$G$ coincides with the differential Galois group taking over the field of the convergent Laurent series $\mathbb{C}\big(\big\{t^{-1}\big\}\big)$. The latter is generated by the group $\mathbb{G}_m\cong \left\{\left(\begin{smallmatrix} c & 0\\ 0& 1/c\end{smallmatrix}\right) \big| \, c\in \mathbb{C}^*\right\}$, $\gamma$ and the two Stokes matrices. The topological monodromy at $t=\infty$ is conjugated to the topological monodromy at $t=0$. Comparing the traces of these matrices yields $-e_1e_2-2={\rm e}^{2\pi {\rm i} \sqrt{8a}}+{\rm e}^{-2\pi {\rm i} \sqrt{8a}}=2\cos \big(2\pi \sqrt{8a}\big)$.

If $e_1e_2\neq 0$, then $G={\rm SL}_2$. Now $e_1e_2=0$ is equivalent to $\sqrt{8a}-\frac{1}{2}\in \mathbb{Z}$ or $a=\frac{(2n+1)^2}{32}$ for some integer $n\geq 0$. In these cases $G$ is contained in a~Borel subgroup of ${\rm SL}_2$ and the operator $\delta_t^2-\big(8a+\frac{1}{2}\delta t^2\big)$ factors as $(\delta_t-F)(\delta_t+F)$ over the field $\mathbb{C}(t)$. In fact $e_1$, $e_2$ are both zero if $a=\frac{(2n+1)^2}{32}$ and $G$ is generated by $\gamma=\left(\begin{smallmatrix}-1& 0\\ 0& -1\end{smallmatrix}\right)$ and
$\mathbb{G}_m\cong \left\{\left(\begin{smallmatrix} c& 0\\ 0& 1/c\end{smallmatrix}\right) \big |\, c\in \mathbb{C}^*\right\}$.

Example: for $a=\frac{1}{32}$ one finds $F=\frac{1}{2}+\sqrt{\delta/2} t$. In fact, a basis of solutions is given by $t^{-1/2}{\rm e}^{\sqrt{\delta /2}t}$, $t^{-1/2}{\rm e}^{-\sqrt{\delta /2} t}$. More generally, a basis of solutions for $a=\frac{(2n+1)^2}{32}$ and (for convenience) with $\delta =2$ is $\big\{t^{-(2n+1)/2}{\rm e}^{t} \big(t^n+\cdots\big), t^{-(2n+1)/2}{\rm e}^{-t}\big(t^n+\cdots\big)\big\}$. The second polynomial is obtained from the first one by changing the sign of the terms~$t^k$ with
$k\equiv (n-1)\bmod 2$.

The second variational equation is
\begin{gather*}
v''+t^{-1}v'-\big(8at^{-2}+\delta/2\big)v=0,\\
 w''+t^{-1}w'-\big(8at^{-2}+\delta/2\big)w=\frac{3}{2}vv''-(v')^2+\frac{3}{2}t^{-1}vv'-\big(16at^{-2}+\delta \big)v^2.
\end{gather*}
For the case $a=\frac{1}{32}$ a MAPLE computation shows that the differential Galois group $G$ of the above two equations has the properties: $G/G^o=C_2$, $G^o/H=\mathbb{G}_m$, $H\cong \mathbb{G}_a^3$ and the action of~$\mathbb{G}_m$ (by conjugation) on~$H$ is not trivial. The details are similar to those of Section~\ref{4.3}. In particular $G^o$ is not commutative. A~similar result holds for all cases $a=\frac{(2n+1)^2}{32}$.

We note that Stoyanova's results in \cite{St2} agree with our's in Section~\ref{4.7} (up to B\"{a}cklund transformations).

\subsection[$P_5\big(\frac{s^2}{2},-\frac{1}{2},-s,-\frac{1}{2}\big)$ and $y=-\frac{t}{s}+1$]{$\boldsymbol{P_5\big(\frac{s^2}{2},-\frac{1}{2},-s,-\frac{1}{2}\big)}$ and $\boldsymbol{y=-\frac{t}{s}+1}$} \label{section4.8}
We note that in \cite[Theorem~2.1(8)]{OO2} there is a typo. According to the ``Clarkson lectures'' \cite[slides~52--53]{C2}, the above choice of parameters corresponds to a Riccati family of solution~$w$ with equation $w'=\frac{s}{t}w^2+\big(1+\frac{1-s}{t}\big)w-\frac{1}{t}$. This Riccati equation has the rational solution $w=1-\frac{t}{s}$.
We refer to \cite[Section~7, especially Section~7.1]{C3} for a derivation of this.

Clarkson's papers \cite{C} and \cite[Section~5.6]{C3} contain a list of rational solutions of $P_5$, e.g., $P_5\big(\frac{1}{2},-\frac{s^2}{2},2-s,-\frac{1}{2}\big)$ with $y=t+s$. We skip these examples since they are equivalent via B\"acklund transformations to one of the two equations in \cite[Theorem~2.1]{OO2}.

The solution $y=-t/s+1$ to $P_5\big(s^2/2,-1/2,-s,-1/2\big)$ induces the VEP $v''+\frac{s-2t}{t(t-s)}v'+\frac{(s-t)^3-s+2t}{t^2(t-s)}v=0$. This equation has a basis of solutions ${\rm e}^{-t}t^{s+1}$, ${\rm e}^tt^{1-s}$ generating the Picard--Vessiot field $\mathbb{C}\big(t,{\rm e}^t,t^{s}\big)$ over $\mathbb{C}(t)$. The corresponding differential Galois group~$G$ is an infinite subgroup of $\mathbb{G}_m\times\mathbb{G}_m$. As in Sections~\ref{4.2} and~\ref{4.3} one computes the second VEP. The resul\-ting inhomogeneous equation yields an extension~$E$ of~$G$ by copies of~$\mathbb{G}_a$, and its connected component~$E^{o}$ is non-commutative.

\subsection[$\deg P_5$ with $\theta_0=\frac{1}{2}$ and solution $y(t)=1-\frac{\theta_1}{2t}$]{$\boldsymbol{\deg P_5}$ with $\boldsymbol{\theta_0=\frac{1}{2}}$ and solution $\boldsymbol{y(t)=1-\frac{\theta_1}{2t}}$} \label{section4.9}

In \cite{OO2} the degenerate fifth Painlev\'e equation $P_5\big(\frac{h^2}{2},-\frac{1}{8},-2,0\big)$ and solution $y=1+\frac{2t^{1/2}}{h}$ is considered (note a small typo in \cite[Theorem~2.1(9)]{OO2}: their `$-8$' should read $-\frac{1}{8}$). Furthermore, \cite[Theorem~2.1(9)]{OO2} together with the last lines of \cite[Section~2.2]{OO2} explain
the relation between~$P_3'(D_6)$ and~$\deg P_5$ (and hence between our Section~\ref{4.4} and the present one). As explained in \cite[p.~9]{APT} the above special~$P_5$ and algebraic solution translate into $\deg P_5\big(\frac{1}{2},\theta_1\big)$ with $\theta_1=h$, and algebraic solution $y=1-\theta_1/(2t)$. It induces as VEP the equation
\[ v''+av'+bv=0 \qquad \text{with} \qquad a=\frac{4t-3\theta_1}{t(2t-\theta_1)},\qquad b= -\frac{32t^3-32t^2\theta_1+8t\theta_1^2+\theta_1}{t^2(2t-\theta_1)}.\]
We note that $h=0$ corresponds to $\theta_1=0$. The given solution has no immediate meaning for $h=0$ and, likewise, the solution $y(t)=1$ has no immediate meaning for the case $\theta_1=0$. Therefore we will suppose that $\theta_1\neq 0$. As pointed out by a referee, $y(t)=1$ has a meaning for an associated Hamiltonian system.

The VEP has three singular points $0$, $\frac{\theta_1}{2}$, $\infty$. The first two singularities are regular singular and $\infty$ is an irregular singularity with Katz invariant~1.

Since $a$ equals $\frac{f'}{f}$ for some $f\in \mathbb{C}(t)$, the differential Galois group $G$ is a subgroup of ${\rm SL}_2$. A standard computation (either by hand or using
MAPLE's DEtools package) shows that at $t=0$, the function $\log t$ is present, implying that~$\mathbb{G}_a\subset G$. The singularity $t=\frac{\theta_1}{2}$ is apparent and ${\rm e}^{4t}$ is present in the formal solutions at $t=\infty$. Thus $G$ also contains a copy of~$\mathbb{G}_m$. If $G\neq {\rm SL}_2$, then $G$ is a~Borel subgroup and the operator $\big(\frac{{\rm d}}{{\rm d}t}\big)^2+a\big(\frac{{\rm d}}{{\rm d}t}\big)+b$ factors (or equivalently the induced Riccati equation has a rational solution). A standard computation shows that in the present case the operator does not factor. One concludes that $G={\rm SL}_2$.

\subsection*{Acknowledgements}

We thank the referees of an earlier version of this paper for their useful suggestions. The first named author thanks the Universidad Simon Bolivar and the Bernoulli Institute of Groningen University for the financial support of his research visit during which the initial version of this paper was written.

\pdfbookmark[1]{References}{ref}
\LastPageEnding


\begin{thebibliography}{99}
\footnotesize\itemsep=0pt

\bibitem{A}
Acosta-Hum\'{a}nez P.B., Nonautonomous {H}amiltonian systems and
 {M}orales-{R}amis theory. {I}.~{T}he case {$\ddot x=f(x,t)$}, \href{https://doi.org/10.1137/080730329}{\textit{SIAM~J.
 Appl. Dyn. Syst.}} \textbf{8} (2009), 279--297, \href{https://arxiv.org/abs/0808.3028}{arXiv:0808.3028}.

\bibitem{APT}
Acosta-Hum\'{a}nez P.B., van~der Put M., Top J., Isomonodromy for the
 degenerate fifth {P}ainlev\'e equation, \href{https://doi.org/10.3842/SIGMA.2017.029}{\textit{SIGMA}} \textbf{13} (2017),
 029, 14~pages, \href{https://arxiv.org/abs/1612.03674}{arXiv:1612.03674}.

\bibitem{CW}
Casale G., Weil J.A., Galoisian methods for testing irreducibility of order two
 nonlinear differential equations, \href{https://doi.org/10.2140/pjm.2018.297.299}{\textit{Pacific~J. Math.}} \textbf{297}
 (2018), 299--337, \href{https://arxiv.org/abs/1504.08134}{arXiv:1504.08134}.

\bibitem{C2}
Clarkson P.A., Painlev\'e equations~-- nonlinear special functions, slides
 presented during the IMA Summer Program Special Functions in the Digital Age,
 Minneapolis, July~22 -- August~2, 2002, available at
 \url{http://www.math.rug.nl/~top/Clarkson.pdf}.

\bibitem{C}
Clarkson P.A., Special polynomials associated with rational solutions of the
 fifth {P}ainlev\'e equation, \href{https://doi.org/10.1016/j.cam.2004.04.015}{\textit{J.~Comput. Appl. Math.}} \textbf{178}
 (2005), 111--129.

\bibitem{C3}
Clarkson P.A., Painlev\'e equations~-- nonlinear special functions, in
 Orthogonal Polynomials and Special Functions, \textit{Lecture Notes in
 Math.}, Vol.~1883, Editors F.~Marcell\'an, W.~Van~Assche, \href{https://doi.org/10.1007/978-3-540-36716-1_7}{Springer}, Berlin,
 2006, 331--411.

\bibitem{G}
Gromak V.I., Laine I., Shimomura S., Painlev\'e differential equations in the
 complex plane, \textit{De Gruyter Studies in Mathematics}, Vol.~28, \href{https://doi.org/10.1515/9783110198096}{Walter de
 Gruyter \& Co.}, Berlin, 2002.

\bibitem{HS}
Horozov E., Stoyanova T., Non-integrability of some {P}ainlev\'e {VI}-equations
 and dilogarithms, \href{https://doi.org/10.1134/S1560354707060056}{\textit{Regul. Chaotic Dyn.}} \textbf{12} (2007), 622--629.

\bibitem{JM}
Jimbo M., Miwa T., Monodromy preserving deformation of linear ordinary
 differential equations with rational coefficients.~{II}, \href{https://doi.org/10.1016/0167-2789(81)90021-X}{\textit{Phys.~D}}
 \textbf{2} (1981), 407--448.

\bibitem{Luk67}
Lukashevich N.A., On the theory of {P}ainlev\'e's third equation,
 \textit{Differ. Uravn.} \textbf{3} (1967), 1913--1923.

\bibitem{Luk68}
Lukashevich N.A., The solutions of {P}ainlev\'e's fifth equation,
 \textit{Differ. Uravn.} \textbf{4} (1968), 1413--1420.

\bibitem{M}
Matsuda M., First-order algebraic differential equations. A~differential
 algebraic approach, \textit{Lecture Notes in Math.}, Vol.~804, \href{https://doi.org/10.1007/BFb0091495}{Springer}, Berlin, 1980.

\bibitem{Mo2}
Morales-Ruiz J.J., A remark about the {P}ainlev\'e transcendents, in Th\'eories
 asymptotiques et \'equations de {P}ainlev\'e, \textit{S\'emin. Congr.},
 Vol.~14, Soc. Math. France, Paris, 2006, 229--235.

\bibitem{MR}
Morales-Ruiz J.J., Ramis J.P., Galoisian obstructions to integrability of
 {H}amiltonian systems, \href{https://doi.org/10.4310/MAA.2001.v8.n1.a3}{\textit{Methods Appl. Anal.}} \textbf{8} (2001),
 33--96.

\bibitem{MRS}
Morales-Ruiz J.J., Ramis J.P., Simo C., Integrability of {H}amiltonian systems
 and differential {G}alois groups of higher variational equations,
 \href{https://doi.org/10.1016/j.ansens.2007.09.002}{\textit{Ann. Sci. \'{E}cole Norm. Sup.~(4)}} \textbf{40} (2007), 845--884.

\bibitem{Mun-vdP}
Muntingh G., van~der Put M., Order one equations with the {P}ainlev\'e
 property, \href{https://doi.org/10.1016/S0019-3577(07)80009-7}{\textit{Indag. Math. (N.S.)}} \textbf{18} (2007), 83--95,
 \href{https://arxiv.org/abs/1202.4633}{arXiv:1202.4633}.

\bibitem{NP}
Nagloo J., Pillay A., On algebraic relations between solutions of a generic
 {P}ainlev\'e equation, \href{https://doi.org/10.1515/crelle-2014-0082}{\textit{J.~Reine Angew. Math.}} \textbf{726} (2017),
 1--27, \href{https://arxiv.org/abs/1112.2916}{arXiv:1112.2916}.

\bibitem{NNPT}
Ngo~Chau L.X., Nguyen K.A., van~der Put M., Top J., Equivalence of differential
 equations of order one, \href{https://doi.org/10.1016/j.jsc.2014.09.041}{\textit{J.~Symbolic Comput.}} \textbf{71} (2015),
 47--59, \href{https://arxiv.org/abs/1303.4960}{arXiv:1303.4960}.

\bibitem{OKSO}
Ohyama Y., Kawamuko H., Sakai H., Okamoto K., Studies on the {P}ainlev\'e
 equations. {V}.~{T}hird {P}ainlev\'e equations of special type {$P_{\rm
 III}(D_7)$} and {$P_{\rm III}(D_8)$}, \textit{J.~Math. Sci. Univ. Tokyo}
 \textbf{13} (2006), 145--204.

\bibitem{OO2}
Ohyama Y., Okumura S., R.~{F}uchs' problem of the {P}ainlev\'e equations from
 the first to the fifth, in Algebraic and Geometric Aspects of Integrable
 Systems and Random Matrices, \textit{Contemp. Math.}, Vol.~593, \href{https://doi.org/10.1090/conm/593/11876}{Amer. Math.
 Soc.}, Providence, RI, 2013, 163--178, \href{https://arxiv.org/abs/math.CA/0512243}{arXiv:math.CA/0512243}.

\bibitem{St3}
Stoyanova T., Non-integrability of {P}ainlev\'e {VI} equations in the
 {L}iouville sense, \href{https://doi.org/10.1088/0951-7715/22/9/008}{\textit{Nonlinearity}} \textbf{22} (2009), 2201--2230.

\bibitem{St2}
Stoyanova T., Non-integrability of {P}ainlev\'e~{V} equations in the
 {L}iouville sense and {S}tokes phenomenon, \href{https://doi.org/10.4236/apm.2011.14031}{\textit{Adv. Pure Math.}}
 \textbf{1} (2011), 170--183.

\bibitem{St1}
Stoyanova T., A note on the {R}.~{F}uchs's problem for the {P}ainlev\'e
 equations, \href{https://arxiv.org/abs/1204.0157}{arXiv:1204.0157}.

\bibitem{St4}
Stoyanova T., Non-integrability of the fourth {P}ainlev\'e equation in the
 {L}iouville--{A}rnold sense, \href{https://doi.org/10.1088/0951-7715/27/5/1029}{\textit{Nonlinearity}} \textbf{27} (2014),
 1029--1044.

\bibitem{SC}
Stoyanova T., Christov O., Non-integrability of the second {P}ainlev\'e
 equation as a {H}amiltonian system, \textit{C.~R.~Acad. Bulgare Sci.}
 \textbf{60} (2007), 13--18, \href{https://arxiv.org/abs/1103.2443}{arXiv:1103.2443}.

\bibitem{U3}
Umemura H., On the irreducibility of the first differential equation of
 {P}ainlev\'e, in Algebraic Geometry and Commutative Algebra, {V}ol.~{II},
 Kinokuniya, Tokyo, 1988, 771--789.

\bibitem{U1}
Umemura H., Second proof of the irreducibility of the first differential
 equation of {P}ainlev\'e, \href{https://doi.org/10.1017/S0027763000001835}{\textit{Nagoya Math.~J.}} \textbf{117} (1990),
 125--171.

\bibitem{U2}
Umemura H., Birational automorphism groups and differential equations,
 \href{https://doi.org/10.1017/S0027763000003111}{\textit{Nagoya Math.~J.}} \textbf{119} (1990), 1--80.

\bibitem{UW}
Umemura H., Watanabe H., Solutions of the second and fourth {P}ainlev\'e
 equations.~{I}, \href{https://doi.org/10.1017/S0027763000006486}{\textit{Nagoya Math.~J.}} \textbf{148} (1997), 151--198.

\bibitem{vdP-Sa}
van~der Put M., Saito M.H., Moduli spaces for linear differential equations and
 the {P}ainlev\'e equations, \href{https://doi.org/10.5802/aif.2502}{\textit{Ann. Inst. Fourier (Grenoble)}}
 \textbf{59} (2009), 2611--2667, \href{https://arxiv.org/abs/0902.1702}{arXiv:0902.1702}.

\bibitem{vdP-Si}
van~der Put M., Singer M.F., Galois theory of linear differential equations,
 \textit{Grundlehren der Mathematischen Wissenschaften}, Vol.~328,
 \href{https://doi.org/10.1007/978-3-642-55750-7}{Springer-Verlag}, Berlin, 2003.

\bibitem{ZF}
\.{Z}o{\l}\c{a}dek H., Filipuk G., Painlev\'e equations, elliptic integrals and
 elementary functions, \href{https://doi.org/10.1016/j.jde.2014.10.018}{\textit{J.~Differential Equations}} \textbf{258} (2015),
 1303--1355.

\end{thebibliography}
\end{document}